\documentclass[report]{owrart}

\usepackage{array}
\usepackage{graphicx}
\graphicspath{{./Figures/}}

\def\vol{\text{vol}}
\def\area{\text{area}}

\def\RR{\mathbb{R}}

\def\bfv{\boldsymbol{v}}
\def\bfx{\boldsymbol{x}}
\def\bfK{\boldsymbol{K}}

\def\bfone{\boldsymbol{1}}

\def\BBB{\mathcal{B}}
\def\SSS{\mathcal{S}}
\def\TTT{\mathcal{T}}
\def\VVV{\mathcal{V}}
\def\EEE{\mathcal{E}}

\newcommand{\red}[1]{{\color{red}#1}}
\newcommand{\blue}[1]{{\color{blue}#1}}
\newcommand{\green}[1]{{\color{green}#1}}
\newcommand{\orange}[1]{{\color{orange}#1}}


\newcommand{\SimS}[1]{\raisebox{-0.75em}{\includegraphics[scale=0.24]{SimplexSplineSmall#1.pdf}}}
\newcommand{\SimSgeneric}{\raisebox{-0.75em}{\includegraphics[scale=0.27]{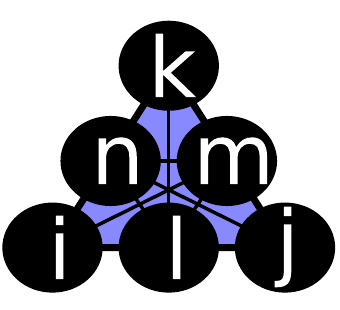}}}

\usepackage{tikz}

\def\PSsmall{
\begin{tikzpicture}[scale = 1.3]
\draw[line width = 0.5] (-0.1, 0.0) -- (0.1,0.0);
\draw[line width = 0.5] (-0.1, 0.0) -- (0.0,0.1732);
\draw[line width = 0.5] ( 0.1, 0.0) -- (0.0,0.1732);
\end{tikzpicture}
}

\def\PSB{
\begin{tikzpicture}[scale = 1.7]
\draw[line width = 0.5] (-0.1, 0.0) -- (0.1,0.0);
\draw[line width = 0.5] (-0.1, 0.0) -- (0.0,0.1732);
\draw[line width = 0.5] ( 0.1, 0.0) -- (0.0,0.1732);

\draw[very thin]  (-0.05, 0.5*0.1732) -- (0.05, 0.5*0.1732);
\draw[very thin]  (-0.05, 0.5*0.1732) -- (0.0, 0.0);
\draw[very thin]  (0.0, 0.0) -- (0.05, 0.5*0.1732);

\draw[very thin]  (0.0, 0.0) -- (0.0,0.1732);
\draw[very thin] (-0.1, 0.0) -- (0.05, 0.5*0.1732);
\draw[very thin] (-0.05, 0.5*0.1732) -- (0.1,0.0);
\end{tikzpicture}
}

\makeatletter
\newenvironment{my_abstract}{%
    \small
    \quotation{\bfseries Abstract.}}
    {\endquotation}
    
\makeatother

\begin{document}

\begin{talk}[Elaine Cohen, Richard Riesenfeld]{Tom Lyche and Georg Muntingh}
{Simplex Spline Bases on the Powell-Sabin 12-Split: Part I}
{Lyche, Tom}

\begin{my_abstract}
We review the construction and a few properties of the S-basis, a simplex spline basis for the $C^1$ quadratic splines on the Powell-Sabin 12-split.
\end{my_abstract}

\bigskip

\noindent Piecewise polynomials or splines defined over triangulations form an indispensable
tool in the sciences, with applications ranging from scattered data fitting  to  finding numerical solutions to partial differential equations. In applications like geometric modelling and solving PDEs by isogeometric methods one often desires a low degree spline with $C^1$, $C^2$ or $C^3$ smoothness. For a general triangulation, it is known that the minimal degree of a triangular $C^r$ element is $4r+1$, e.g., degrees
5, 9, 13 for the classes $C^1$, $C^2$ or $C^3$. To obtain smooth splines of lower degree
one can split each triangle in the triangulation into several subtriangles. One
such split that we consider here is the Powell-Sabin 12-split of a triangle.
\begin{center}
\includegraphics[scale=0.6]{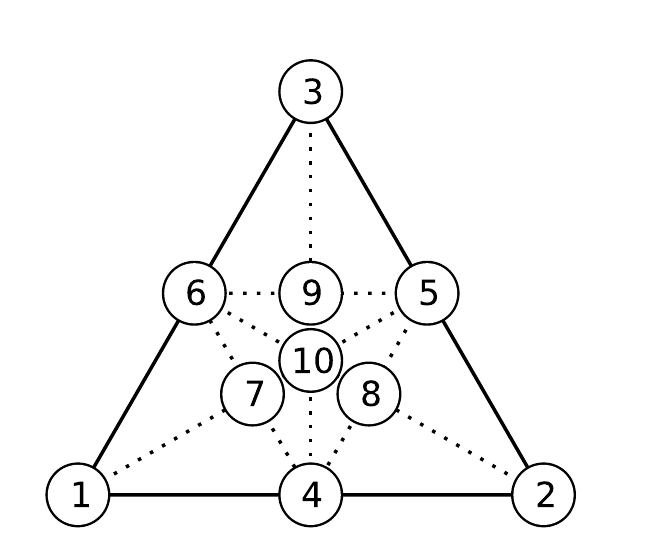}\\The 12-split with numbering of vertices.
\end{center}

Once a space is chosen one determines its dimension. The spaces $\SSS^1_2(\PSB)$ and $\SSS^3_5(\PSB)$ of $C^1$ quadratics and $C^3$ quintics on the 12-split $\PSB$ of a single triangle have dimension 12 and 39, respectively. Over a general triangulation $\TTT$ of a polygonal domain we can 12-split each triangle in $\TTT$ to obtain a triangulation $\TTT_{12}$. The  dimensions of the corresponding $C^1$ quadratic and $C^2$ quintic spaces (the latter with $C^3$ supersmoothness  at the vertices and the interior edges of each macro triangle) are $3|\VVV| + |\EEE|$ and $10|\VVV| + 3|\EEE|$, respectively, where $|\VVV|$ and $|\EEE|$ are the number of vertices and edges in $\TTT$. Moreover, in addition to giving $C^1$ and $C^2$ spaces on any triangulation these spaces are suitable for multiresolution analysis, see for example \cite{TL_LycheMuntingh14}.

To compute with these spaces one needs a suitable basis. In the univariate case the B-spline basis is an obvious choice. In this talk we consider a bivariate generalization known as simplex splines. We review the construction and a few  properties shown in \cite{TL_CohenLycheRiesenfeld13} of the S-basis consisting of $C^1$ quadratic simplex splines in $\SSS^1_2(\PSB)$. We also introduce some concepts needed in Part II of this talk given by Georg Muntingh. 

{\bf A short background on simplex splines}.
Let $\bfK = \{\bfv_1\cdots \bfv_{d+s+1}\} \subset \RR^s$ be a finite multiset. Consider a simplex $\sigma = [\overline{\bfv}_1, \ldots, \overline{\bfv}_{d+s+1}] \subset \RR^{d+s}$ together with a projection $\pi: \sigma \longrightarrow \RR^s$ satisfying $\pi(\overline{\bfv}_i) = \bfv_i$. We define the \emph{simplex spline}
$ B[\bfK](\bfx) = \vol_d\big(\sigma\cap \pi^{-1}(\bfx)\big)/\vol_{d+s} (\sigma)$. For instance, three knots in $\RR^1$ define a linear B-spline, four knots in $\RR^1$ define a quadratic B-spline, and four knots in $\RR^2$ define a linear bivariate simplex spline: 
\begin{center}
\includegraphics[scale=0.7, clip = true, trim = 0 10 25 16]{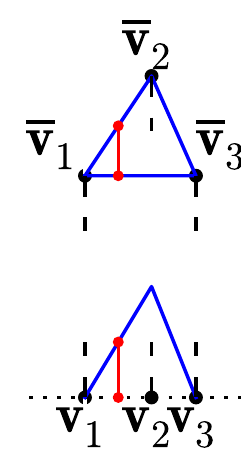}\qquad
\includegraphics[scale=0.7, clip = true, trim = 0 10 25 22]{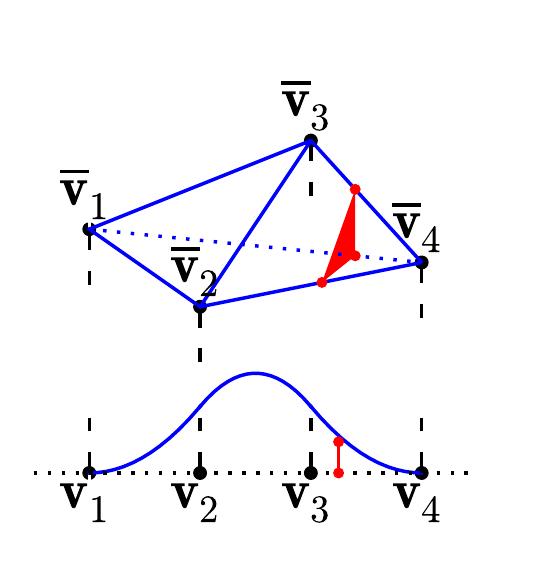}\qquad
\includegraphics[scale=0.7, clip = true, trim = 0 10 25 21]{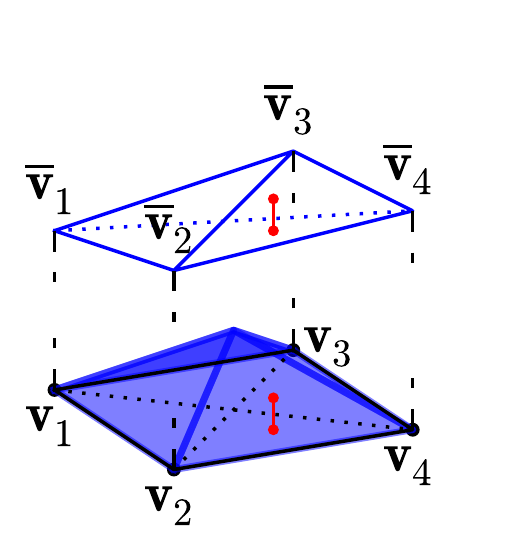}
\end{center}

Simplex splines have all the usual properties of univariate B-splines. This includes continuity which can be controlled locally, a recurrence relation, and differentiation and knot insertion formulas. The support of a simplex spline is the convex hull of its knots, and in $\RR^2$ the collection of knotlines is obtained by connecting  each knot to all other knots (the complete graph). A simplex spline with $d+3$ knots in $\RR^2$ has $d-m+1$ continuous derivatives across a knot line containing $m$ knots counting multiplicites. 

{\bf Simplex splines on the 12-split}.
Since the knotlines form a complete graph the simplex splines are natural candidates for a $C^r$ basis on this split. A simplex spline on the 12-split will have a knotset of the form $\bfK = \{\bfv_1^{m_1} \cdots \bfv_{10}^{m_{10}}\}$, where $\bfv_1,\ldots,\bfv_{10}$ are the vertices numbered as above, and $m_i\ge 0$ is the multiplicity of $\bfv_i$, i.e., the number of repetitions of $\bfv_i$ in the multiset. A convenient scaling is the 
\emph{(area normalized) simplex spline} $Q[\bfK]: \RR^2 \longrightarrow \RR$, recursively defined by
\[
Q[\bfK](\bfx) :=
\left\{ \begin{array}{cl}
0 & \text{if~}\area([\bfK]) = 0,\\
\bfone_{[\bfK)}(\bfx)\frac{\area(\PSsmall)}{\area([\bfK])} & \text{if~}\area([\bfK]) \neq 0\text{~and~} |\bfK| = 3,\\
\sum_{j = 1}^{10} \beta_j Q[\bfK\backslash \bfv_j](\bfx) &  \text{if~}\area([\bfK]) \neq 0\text{~and~} |\bfK| > 3,\\
\end{array} \right. \]
with $\bfx = \beta_1 \bfv_1 + \cdots + \beta_{10} \bfv_{10}, \beta_1 + \cdots + \beta_{10} = 1$, and $\beta_i = 0$ whenever $m_i = 0$.

By Theorem 4 in \cite{TL_Micchelli79} this definition is independent of the choice of the $\beta_j$.  
Whenever $m_7 = m_8 = m_9 = m_{10} = 0$, 
we use the graphical notation
\[ \SimSgeneric := Q[\bfv_1^i \bfv_2^j \bfv_3^k \bfv_4^l \bfv_5^m \bfv_6^n]. \]

{\bf B-splines on the boundary}. It is useful for the simplex splines to restrict to consecutive univariate B-splines on the boundary. For example, on  $[\bfv_1,\bfv_2]$  the quadratic simplex splines 
$\frac14$\SimS{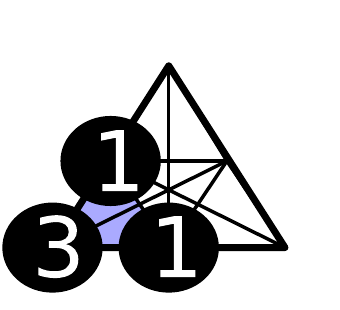}, $\frac12$\SimS{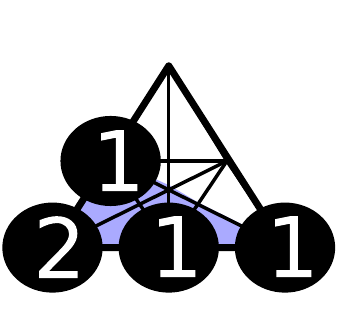},$\frac12$\SimS{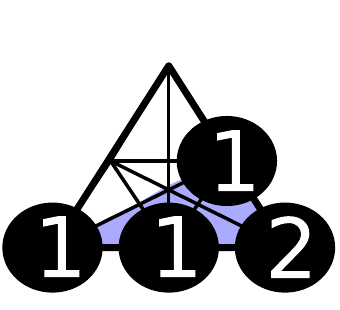}, $\frac14$\SimS{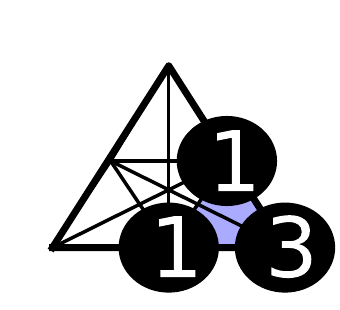} restrict to:
\vskip 0.1cm

\begin{tabular}{m{11em}c}
\includegraphics[scale=0.3]{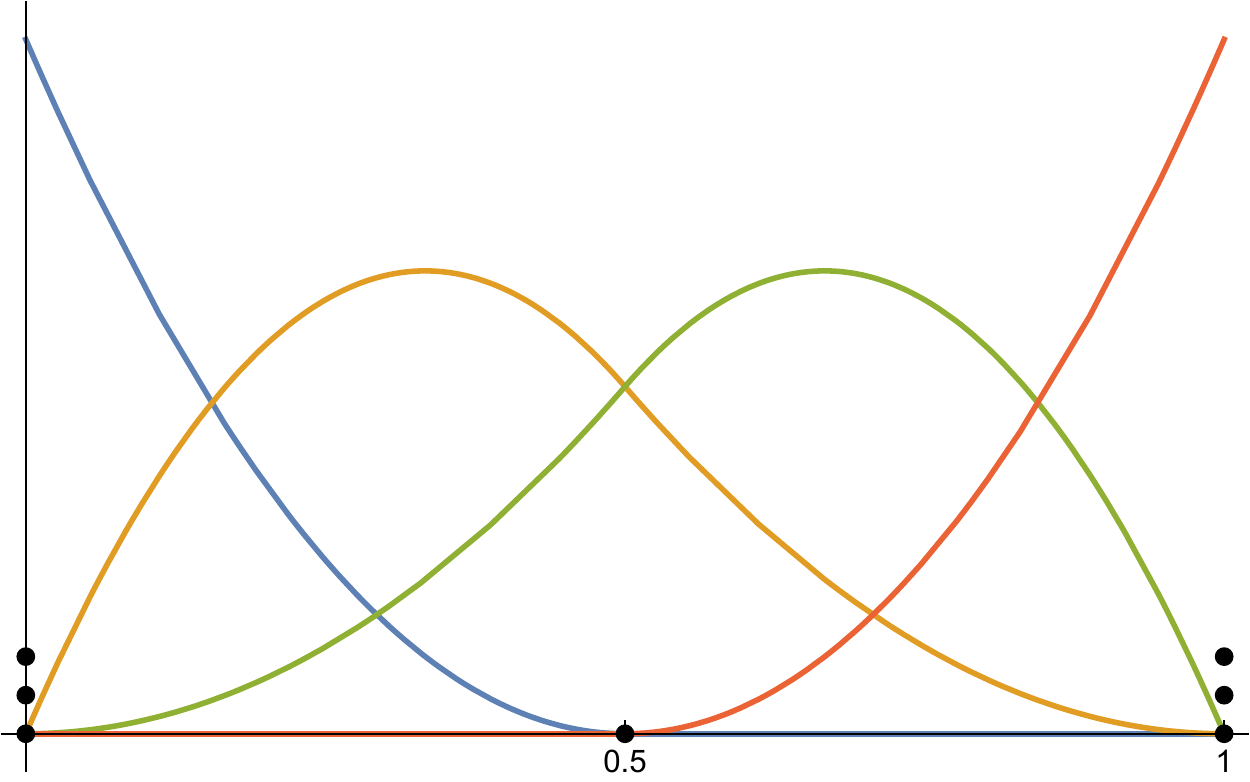}&$\blue{B[\bfv_1^3,\bfv_4]}$, $\orange{B[\bfv_1^2,\bfv_4,\bfv_2]}$, $\green{B[\bfv_1,\bfv_4,\bfv_2^2]}$, $\red{B[\bfv_4,\bfv_2^3]}$
\end{tabular}
\vskip 0.3cm

{\bf Symmetries}.
Identifying a triangle  with an equilateral triangle, its symmetries
\begin{center}
\includegraphics[scale=0.55]{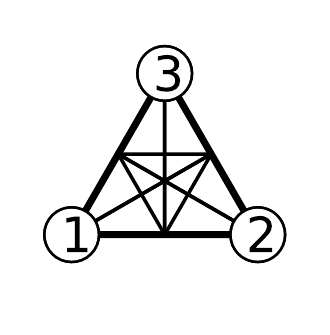}
\includegraphics[scale=0.55]{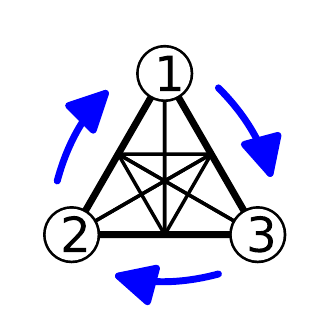}
\includegraphics[scale=0.55]{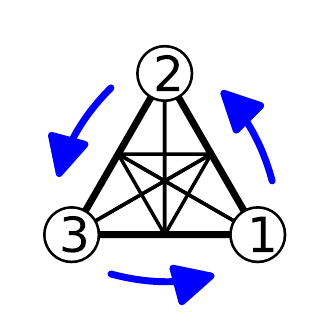}
\includegraphics[scale=0.55]{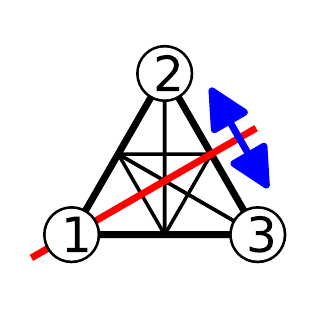}
\includegraphics[scale=0.55]{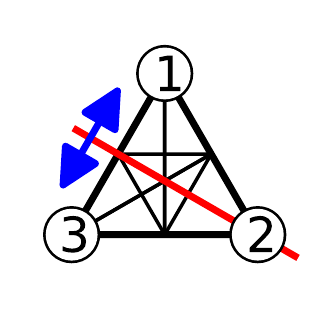}
\includegraphics[scale=0.55]{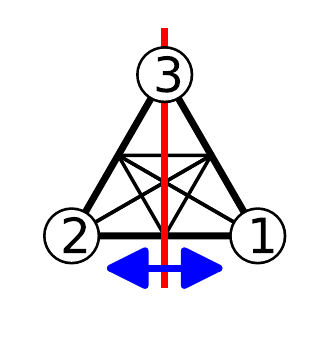}
\end{center}
form a group $S_3$ that acts on the simplex splines by permuting knots. We write
\[ [\BBB]_{S_3} := \{Q[\sigma(\bfK)]\,:\, Q[\bfK]\in \BBB,\ \sigma\in S_3\} \]
for the set of simplex splines related to $\BBB$ by a symmetry in $S_3$.
Let
\begin{equation*}
\begin{aligned}
& c_4 := \frac{c_1 + c_2}{2},\quad c_5 := \frac{c_2 + c_3}{2},\quad c_6 := \frac{c_1 + c_3}{2},\\
& c_7 := \frac{c_4 + c_6}{2},\quad c_8 := \frac{c_4 + c_5}{2},\quad c_9 := \frac{c_5 + c_6}{2},\quad
c_{10} := \frac{c_1 + c_2 + c_3}{3}.
\end{aligned}
\end{equation*}
Via the identification $c_i \leftrightarrow \bfv_i$ with the vertices of $\PSB$, the group $S_3$ acts on polynomials in $c_1,\ldots,c_{10}$ and simplex splines, or combinations of these, e.g.,
\begin{align*}
\left[c_4c_{10}\SimS{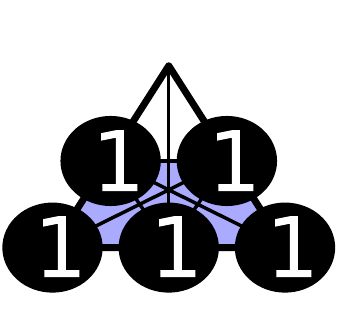}\right]_{S_3} & = \left\{c_4c_{10}\SimS{110111}, c_5c_{10}\SimS{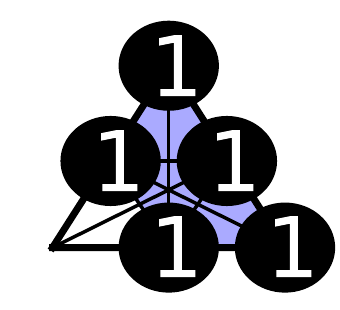}, c_6c_{10}\SimS{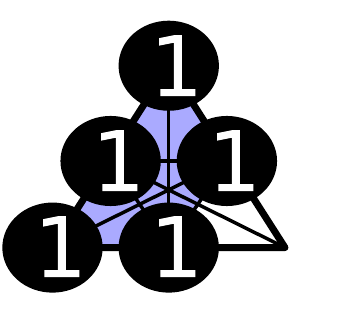}\right\}.
\end{align*}
{\bf The quadratic S-basis}. It is given by 
\[\left[\frac14\SimS{300101}, \frac12\SimS{210101}, \frac34\SimS{110111}\right]_{S_3}
= \left\{\frac14\SimS{300101}, \frac14\SimS{030110}, \frac14\SimS{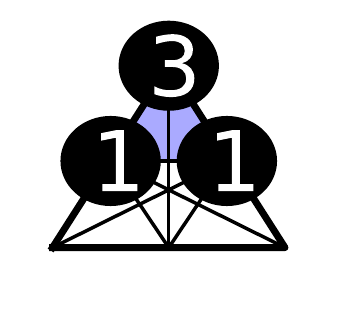},\ldots,
\frac34\SimS{011111}, \frac34\SimS{110111}
\right\}
\]
and is the unique simplex spline basis for $\SSS_2^1(\PSB)$ with local linear independence. Moreover, it is symmetric, reduces to B-splines on the boundary, can be computed by a pyramidal scheme, and has B\'ezier-like smoothness conditions across adjacent macro triangles. Furthermore, it has a barycentric Marsden identity 
\[
 \left( \sum \left[c_1\SimS{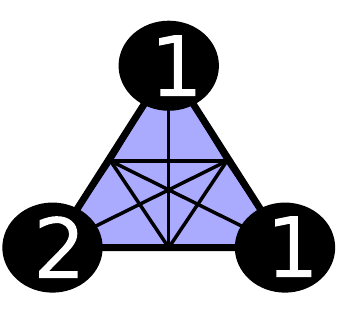}\right]_{S_3}\right)^2 =
 \sum \left[\frac14 c_1^2\SimS{300101}\right]_{S_3} \cup \left[\frac34 c_4c_{10} \SimS{110111}\right]_{S_3} \cup \left[\frac12 c_1 c_4 \SimS{210101} \right]_{S_3},
\]
which yields polynomial reproduction, explicit dual functionals and a simple quasi-interpolant. These show that the S-basis is stable independently of the geometry, which implies an $h^2$ bound on the distance between a spline and its control surface.

\end{talk}

\begin{thebibliography}{99}
\bibitem{TL_CohenLycheRiesenfeld13}
Elaine Cohen, Tom Lyche, Richard Riesenfeld, \textit{A B-spline-like basis for the Powell-Sabin 12-split based on simplex splines}, Math. Comp. \textbf{82} (2013), no. 283, 1667--1707.

\bibitem{TL_LycheMuntingh14}
Tom Lyche and Georg Muntingh, \textit{A Hermite interpolatory subdivision scheme for $C^2$-quintics on the Powell-Sabin 12-split}, Comput. Aided Geom. Design \textbf{31} (2014), no. 7--8, 464--474.

\bibitem{TL_LycheMuntingh15}
Tom Lyche and Georg Muntingh, \textit{Stable simplex spline bases for $C^3$ quintics on the Powell-Sabin 12-split}, Available at \texttt{http://arxiv.org/abs/1504.02628}.

\bibitem {TL_Micchelli79}
Charles A. Micchelli, 
   \textit{On a numerically efficient method for computing multivariate
   $B$-splines}, in 
   "Multivariate approximation theory", Walter Schempp and Karl Zeller (eds.), International Series of  Numerical Mathematics Vol. 51, Birkh\"auser Verlag, Basel, Boston, Stuttgart,
1979, 211--248.
\end{thebibliography}
\end{document}